\documentclass[12pt]{article}
\usepackage{amssymb}
\usepackage{amsmath}
\usepackage{amstext}
\usepackage{amsfonts}
\usepackage{amsthm}
\usepackage{amscd}
\numberwithin{equation}{section}

\swapnumbers
\theoremstyle{plain}
\newtheorem{thm}{Theorem}[section]
\newtheorem{lem}[thm]{Lemma}
\newtheorem{prop}[thm]{Proposition}
\newtheorem{cor}[thm]{Corollary}
\newtheorem{definition}[thm]{Definition}

\theoremstyle{remark}
\newtheorem{rem}[thm]{Remark}

\def\cO{\mathcal{O}}

\def\map#1{\ \smash{\mathop{\longrightarrow}\limits^{#1}}\ }

\def\Sing{\mathrm{Sing}}

\def\LL{\mathcal{L}}

\def\MM{\mathbf{M}}
\def\N{\mathcal{N}}

\def\PP{\mathbb{P}}

\def\ZZ{\mathbb{Z}}
\def\TT{\mathbb{T}}

\def\FF{\mathbb{F}}
\def\EE{\mathbb{E}}
\def\AA{\mathbb{A}}
\def\CC{\mathbb{C}}
\def\WW{\mathbb{W}}
\def\DD{\mathcal{D}}

\def\T{\mathcal{T}}
\def\BL{\mathcal{B}l}
\def\H{\mathcal{H}}
\def\Hess{\mathrm{Hess}}
\def\St{\mathrm{St}}
\def\P{\mathcal{P}}
\def\R{\mathcal{R}}

\def\deg{\mathrm{deg} \:}

\def\sym{\mathrm{Sym}}
\def\ker{\mathrm{ker} \:}
\def\pic{\mathrm{Pic}}
\def\mult{\mathrm{mult}}

\def\dim{\mathrm{dim} \:}
\def\det{\mathrm{det} \:}
\def\rk{\mathrm{rk} \:}

\def\im{\mathrm{im} \:}

\def\Hom{\mathcal{H}\mathrm{om}}
\def\supp{\mathrm{supp} \:}

\def\ext{\mathrm{Ext}}
\def\End{\mathrm{End}}

\def\lra{\longrightarrow}

\def\ra{\rightarrow}

\def\lms{\longmapsto}

\def\GL{\mathrm{GL}}
\def\SL{\mathrm{SL}}
\def\Sp{\mathrm{Sp}}

\def\pp{\mathbb{P}}

\def\M{\mathcal{M}}
\def\N{\mathcal{N}}
\def\U{\mathcal{U}}
\def\K{\mathcal{K}}
\def\E{\mathcal{E}}
\def\F{\mathcal{F}}

\title{Self-duality of Coble's quartic hypersurface and applications}
\author{Christian Pauly}
\begin{document}
\maketitle

\begin{abstract}
The moduli space $\M_0$ of semi-stable rank $2$ vector bundles with
fixed trivial determinant over a non-hyperelliptic curve $C$ of genus
$3$ is isomorphic to a quartic hypersurface in $\pp^7$ (Coble's
quartic). We show that $\M_0$ is self-dual and that its polar map
associates to a stable bundle $E \in \M_0$ a bundle $F$ which is
characterized by $\dim H^0(C, E \otimes F) =4$. The
projective space $\pp H^0(C, E \otimes F)$ is equipped with a net
of quadrics $\Pi$ and it is shown that the map which associates 
to $E \in \M_0$ the isomorphism class of the plane quartic
Hessian curve of $\Pi$ is a dominant map to the moduli
space of genus $3$ curves. 
\end{abstract}

\section{Introduction}
In his book \cite{co} A.B. Coble constructs for any non-hyperelliptic
curve $C$ of genus $3$ a quartic hypersurface in $\pp^7$ which 
is singular along the Kummer variety $\K_0 \subset \pp^7$ of 
$C$. It is shown in
\cite{narram} that this hypersurface is isomorphic to the 
moduli space $\M_0$ of semi-stable rank $2$  vector bundles
with fixed trivial determinant. For many reasons Coble's quartic
hypersurface may be viewed as a genus-$3$-analogue 
of a Kummer surface,i.e., a quartic surface $S \subset \pp^3$ with 
$16$ nodes. For example 
the restriction of $\M_0$ to an eigenspace $\pp^3_\alpha \subset
\pp^7$ for the action of a $2$-torsion point $\alpha \in JC[2]$ 
is isomorphic to a Kummer surface (of the corresponding Prym variety). 
It is classically known (see e.g. \cite{gh}) that a 
Kummer surface $S \subset \pp^3$ is self-dual.

\bigskip

In this paper we show that this property also holds for the
Coble quartic $\M_0$ (Theorem \ref{maintheorem}). The rational
polar map $\DD: \pp^7 \lra \left( \pp^7 \right)^*$ maps
$\M_0$ birationally to $\M_\omega \subset \left( \pp^7 \right)^*$,
where $\M_\omega$ ($\cong \M_0$) is the moduli space
parametrizing vector bundles with fixed canonical determinant.
More precisely we show that the embedded tangent space at a
stable bundle $E$ to $\M_0$ corresponds to a semi-stable bundle
$\DD(E) = F \in \M_\omega$, which is characterized by the
condition $\dim H^0(C, E\otimes F) = 4$ (its maximum).
We also show that $\DD$ resolves to a morphism $\widetilde{\DD}$
by two successive blowing-ups, and that $\DD$ contracts the
trisecant scroll of $\K_0$ to the Kummer variety $\K_\omega \subset
\M_\omega$.

\bigskip

The condition which relates $E$ to its ``tangent space bundle'' $F$,
namely $\dim H^0(C, E\otimes F) = 4$, leads to many geometric
properties. First we observe that $\pp H^0(C, E\otimes F)$
is naturally equipped with a net of quadrics $\Pi$ whose base 
points (Cayley octad) correspond bijectively to the $8$
line subbundles of maximal degree of $E$ (and of $F$). The Hessian
curve $\Hess(E)$ of the net of quadrics $\Pi \cong |\omega|^*$
is a plane quartic curve, which is everywhere tangent (Proposition
\ref{evtan}) to the canonical curve $C \subset |\omega|^*$,i.e.,
$\Hess(E) \cap C = 2\Delta(E)$ for some divisor $\Delta(E) \in
|\omega^2|$. Since these constructions are $JC[2]$-invariant,
we introduce the quotient $\N = \M_0/JC[2]$ parametrizing
$\pp \SL_2$-bundles over $C$ and we show (Proposition \ref{proj})
that the map $\N \map{\Delta} |\omega^2|$, $E \lms \Delta(E)$ 
is the restriction of the projection from the 
projective space $\N \subset |\overline{\LL}|^*
= \pp^{13}$ ($\overline{\LL}$ is the ample generator of $\pic(\N)$)
with center of projection given by the linear span of the 
Kummer variety $\K_0 \subset \N$ ($\K_0$ parametrizes 
decomposable $\pp \SL_2$-bundles).

\bigskip

We also show (Corollary \ref{dom}) that the Hessian map $\N \ra \R$, 
$E \mapsto \Hess(E)$ is finite of degree $72$, where
$\R$ is the rational space parametrizing plane quartics everywhere
tangent to $C \subset |\omega|^* = \pp^2$. Considering the
isomorphism class of $\Hess(E)$, we  deduce that the map 
$\Hess: \N \ra \M_3$ is dominant, where $\M_3$ is the moduli
space of smooth genus $3$ curves. We actually prove that
some Galois-covers $\widetilde{\N} \ra \N$ and $\P_C \ra \R$ are
birational (Proposition \ref{birat}). In  particular we endow the
space $\widetilde{\N}$, parametrizing $\pp \SL_2$-bundles $E$
with an ordered set of $8$ line subbundles of $E$ of
maximal degree, with an action of the Weyl group  $W(E_7)$ such that
the action of the central element $w_0 \in W(E_7)$ coincides 
with the polar map $\DD$.

\bigskip

We hope that these results will be useful for dealing with 
several open problems, e.g. rationality of the moduli spaces
$\M_0$ and $\N$.

\bigskip

I would like to thank S. Ramanan for some inspiring 
discussions on Coble's quartic.

\section{The geometry of Coble's quartic}

In this section we briefly recall some known results related to 
Coble's quartic hypersurface, which can be found in the 
literature, e.g. \cite{do}, 
\cite{la2}, \cite{narram}, \cite{opp}. We refer to \cite{b1}, 
\cite{b2} for the results on the geometry of the moduli of
rank $2$ vector bundles. 

\subsection{Coble's quartic as moduli of vector bundles}

Let $C$ be a smooth non-hyperelliptic curve of genus $3$ with 
canonical line bundle $\omega$. Let $\pic^{d}(C)$ be the Picard variety
parametrizing degree $d$ line bundles over $C$ and $JC:= \pic^0(C)$ be
the Jacobian variety. We denote by $\K_0$ the Kummer variety of $JC$ and
by $\K_\omega$ the quotient of $\pic^2(C)$ by the involution $\xi \mapsto
\omega \xi^{-1}$.  Let $\Theta \subset \pic^2(C)$ be the Riemann
Theta divisor and let $\Theta_0 \subset JC$ be a symmetric Theta
divisor, i.e., a translate of $\Theta$ by a theta-characteristic.
We also recall that the two linear systems  $|2\Theta|$
and $|2\Theta_0|$ are canonically dual to each other,i.e., 
$|2\Theta|^* \cong |2\Theta_0|$.

\bigskip
\noindent
Let $\M_0$ (resp. $\M_\omega$) denote the moduli space of semi-stable
rank $2$ vector bundles over $C$ with fixed trivial (resp. canonical)
determinant. The singular locus of $\M_0$ is isomorphic to $\K_0$ and points
in $\K_0$ correspond to bundles $E$ whose $S$-equivalence 
class $[E]$ contains a decomposable bundle of the form 
$M \oplus M^{-1}$ for $M \in JC$. We have natural 
morphisms 
$$ \M_0 \map{D} |2\Theta| = \pp^7, \qquad \M_\omega \map{D} 
|2\Theta_0| = |2\Theta|^*,$$
which send a stable bundle $E \in \M_0$ to the divisor $D(E)$ whose
support equals the set $\{ L \in \pic^2(C) \ | \ \dim H^0(C, E \otimes L) > 0
\}$ (if $E\in \M_\omega$, replace $\pic^2(C)$ by $JC$). On the
semi-stable boundary $\K_0$ (resp. $\K_\omega$) the morphism
$D$ restricts to the Kummer map. The moduli spaces $\M_0$ and $\M_\omega$
are isomorphic, although non-canonically (consider tensor product
with a theta-characteristic). It is known that the Picard group
$\pic(\M_0)$ is $\ZZ$ and that $|\LL|^* = |2\Theta|$, where $\LL$ is the
ample generator of $\pic(\M_0)$.

\bigskip
\noindent
The main theorem of \cite{narram} asserts that $D$ embeds $\M_0$ as
a quartic hypersurface in $|2\Theta|= \pp^7$, which was originally described
by A.B. Coble \cite{co} (section 33(6)). 
Coble's quartic is characterized by a uniqueness
property: it is the unique (Heisenberg-invariant)
quartic which is singular along the Kummer variety $\K_0$ 
(see \cite{la2} Proposition 5).

\bigskip
\noindent
We recall that Coble's quartic hypersurfaces $\M_0 \subset |2\Theta|$
and $\M_\omega \subset |2\Theta_0|$ contain some distinguished points.
First (\cite{co} section 48(4), \cite{la1}, \cite{opp}) there exists
a unique stable bundle $A_0 \in \M_\omega$ such that $\dim H^0(C,A_0)
=3$ (its maximal dimension). We define for any theta-characteristic
$\kappa$ and for any $2$-torsion point $\alpha \in JC[2]$ the stable
bundles, called {\em exceptional} bundles
\begin{equation} \label{distbun}
A_\kappa := A_0 \otimes \kappa^{-1} \in \M_0 \qquad \text{and} \qquad
A_\alpha := A_0 \otimes \alpha \in \M_\omega.
\end{equation}

\subsection{Global and local equations of Coble's quartic}

Let $F_4$ be the Coble quartic,i.e., the equation of $\M_0 \subset
|2\Theta| = \pp^7$. Then the eight partials $C_i = 
\frac{\partial F_4}{\partial X_i}$ for $1 \leq i \leq 8$ ($X_i$ are 
coordinates for $|2\Theta|$) define the Kummer variety $\K_0$
scheme-theoretically (\cite{la2} Theorem IV.6). We also need
the following results (\cite{la2} Theorem 6 bis).
\begin{itemize}
\item[(i)] The \'etale local equation (in affine space $\AA^7$) of 
Coble's quartic at the point
$[\cO \oplus \cO]$ is $T^2 = \det [T_{ij}]$, with coordinates
$T,  T_{ij}$ with $T_{ij} = T_{ji}$ and $1 \leq i,j \leq 3$.
\item[(ii)] The \'etale local equation at the point $[M \oplus M^{-1}]$ with
$M^2 \not= \cO$ is a rank $4$ quadric $\det [T_{ij}] = 0$, where
$T_{ij}$, $1 \leq i,j \leq 2$ are four coordinates on $\AA^7$.
\end{itemize}

Hence any point $[M \oplus M^{-1}] \in \K_0$ has multiplicity $2$
on $\M_0$.

\subsection{Extension spaces}

Given $L \in \pic^1(C)$ we introduce the $3$-dimensional 
space $\pp_0(L) := |\omega L^2|^* = \pp \ext^1(L,L^{-1})$. A point 
$e \in \pp_0(L)$ corresponds to an isomorphism class of extensions
\begin{equation} \label{exte}
0 \lra L^{-1} \lra E \lra L \lra 0 \qquad (e) 
\end{equation}
and the composite of the classifying map $\pp_0(L) \ra \M_0$ followed
by the embedding $D : \M_0 \ra |2\Theta|$ is linear and injective
(\cite{b2} Lemme 3.6). It is shown that a point $e \in \pp_0(L)$ represents
a stable bundle precisely away from $\varphi(C)$, where $\varphi$ is
the map induced by the linear system $|\omega L^2|$. A point $e = \varphi(p)$
for $p \in C$ is represented by the decomposable bundle $L(-p) \oplus
L^{-1}(p)$.

\bigskip
\noindent
We also introduce the projective spaces 
$\pp_\omega(L) := |\omega^2 L^{-2}|^* =
\pp \ext^1(\omega L^{-1},L)$. A point $f \in \pp_\omega(L)$
corresponds to an extension 
\begin{equation} \label{extf}
0 \lra L \lra F \lra \omega L^{-1} \lra 0 \qquad (f)
\end{equation}
Similarly, we have an injective classifying map 
$\pp_\omega(L) \ra \M_\omega$. Although we will not use this
fact, we observe that $\pp_0(L) = \pp_\omega(\kappa L^{-1})$ for
any theta-characteristic $\kappa$.

\bigskip
\noindent

It is well-known that the Kummer variety $\K_0 \subset |2\Theta|$ admits
a $4$-dimensional family of trisecant lines. It follows from \cite{opp}
Theorem 1.4 and Theorem 2.1 that any trisecant line to $\K_0$ is
contained in some space $\pp_0(L)$, where it is a trisecant to 
the curve $\varphi(C) \subset \pp_0(L)$. We denote by $\T_0$ the
trisecant scroll, which is a divisor in $\M_0$. Similarly we 
define $\T_\omega \subset \M_\omega$.

\bigskip
\noindent

The main tool for the proof of the self-duality is that $\M_0$ 
(resp. $\M_\omega$) can be
covered by the projective spaces $\pp_0(L)$ (resp. $\pp_\omega(L)$). 
This is expressed by the following result \cite{narram}(see also \cite{op}): 
there exist a rank $4$ 
vector bundles $\U_0$ and $\U_\omega$ over $\pic^1(C)$ 
such that $\forall L \in 
\pic^1(C)$, $\left(\pp \U_0 \right)_L \cong  \pp_0(L)$,
 $\left(\pp \U_\omega \right)_L \cong  \pp_\omega(L)$ and
their associated classifying morphisms $\psi_0$ and $\psi_\omega$,
$$
\begin{CD}
\pp \U_0 @>\psi_0>> \M_0 \subset |2\Theta| \\
@VVV  \\
\pic^1(C)
\end{CD}
\qquad
\begin{CD}
\pp \U_\omega @>\psi_\omega>> \M_\omega \subset |2\Theta_0| \\
@VVV  \\
\pic^1(C)
\end{CD}
$$
are  surjective (Nagata's theorem) and of degree $8$ (see section 4.1).

\subsection{Tangent spaces to Theta-divisors}

Following \cite{b2} section 2, we associate to any $[F] \in \M_\omega \subset
|2\Theta_0|$ the divisor $\Delta(F) \subset \M_0 \subset |2\Theta|$
which has the properties
\begin{enumerate}
\item $\supp \Delta(F) = \{ [E] \in \M_0 \ | \ \dim 
H^0(C, E \otimes F) > 0 \},$
\item $\Delta(F) \in |\LL| \cong |2\Theta|^*$ is mapped to $[F]$ under
the canonical duality $|2\Theta|^* \cong |2\Theta_0|.$
\end{enumerate}
Symmetrically, we associate to any $E \in \M_0$ the divisor
$\Delta(E) \subset \M_\omega$ with the analoguous properties. 

\bigskip
\noindent
For any $E,F$ with $[E] \in \M_0$ and $[F] \in \M_\omega$, the 
rank $4$ vector bundle $E \otimes F = \Hom(E,F)$ is equipped
with a $\omega$-valued non-degenerate quadratic form (given by
the determinant of local sections), hence, by Mumford's 
parity theorem \cite{mum}, $\dim H^0(C, E\otimes F)$ is even.
The divisor $\Delta(E)$ is defined as the Pfaffian divisor
associated to a family $E \otimes \mathcal{F}$ of orthogonal 
bundles \cite{ls}. Adapting the proof of \cite{la1} Lemme II.2
to Pfaffian divisors, we see that
$$ \mult_{[E]} \Delta(F) \geq \frac{1}{2} \dim H^0(C,E \otimes F).$$
\begin{lem} \label{tansing}
Suppose that $E$ is stable and that $\dim H^0(C, E\otimes F) \geq 4$.
Then $\Delta(F) \subset \M_0$ is singular at $E$ and the
embedded tangent space $\TT_E \M_0 \in |2\Theta|^* \cong |2\Theta_0|$
corresponds to the point $[F] \in |2\Theta_0|$.
\end{lem}

\begin{proof}
The first assertion is an immediate consequence of the previous
inequality. To show the second, it is enough to observe that, since
$E$ is a singular point of the divisor $\Delta(F)$, we have equality
between the Zariski tangent spaces $T_E \Delta(F) = T_E \M_0$ and
$T_E \Delta(F)$ coincides with the hyperplane cutting out the 
divisor $\Delta(F)$, which corresponds to the point $[F]$ by 
property (2).
\end{proof}

We will also need the dual version.
\begin{lem} \label{tansingdual}
Suppose that $F$ is stable and that $\dim H^0(C, E\otimes F) \geq 4$.
Then $\Delta(E) \subset \M_\omega$ is singular at $F$ and the
embedded tangent space $\TT_F \M_\omega \in |2\Theta_0|^* \cong |2\Theta|$
corresponds to the point $[E] \in |2\Theta|$.
\end{lem}

\section{Self-duality}

\subsection{Statement of the main theorem}

Let $\DD$ be the rational map defined by the polars of Coble's
quartic $F_4$,i.e., the eight cubics $C_i$,
\begin{eqnarray*} 
\DD :  |2\Theta| &  \lra &  |2\Theta|^* \cong |2\Theta_0| \\
 \cup \  & & \  \cup \\
 \M_0  & & \M_\omega  
\end{eqnarray*}
Note that $\DD$ is defined away from $\K_0$.
Geometrically, $\DD$ maps a stable bundle $E \in \M_0$ to
the hyperplane defined by the embedded tangent space $\TT_E \M_0$
at the smooth point $E$.
The main theorem of this paper is the following
\begin{thm}[Self-duality] \label{maintheorem}
The moduli space $\M_0$ is birationally mapped by $\DD$ to $\M_\omega$,i.e.,
$\M_\omega$ is the dual hypersurface of $\M_0$. More precisely, we have
\begin{enumerate}
\item 
$\DD$ restricts to an isomorphism $\M_0 \setminus \T_0 \map{\sim}
\M_\omega \setminus \T_\omega$.
\item 
$\DD$ contracts the divisor $\T_0$ ($\in |\LL^8|$) to $\K_\omega$.
\item
For any stable $E \in \M_0$, the moduli point $\DD(E) \in \M_\omega$ can
be represented by a semi-stable bundle $F$, which satisfies 
$\dim H^0(C,E \otimes F) \geq 4$. Moreover, if $E \in \M_0 \setminus
\T_0$ then there exists a unique stable bundle $F = \DD(E)$
for which $\dim H^0(C, E\otimes F)$ has its
maximal value $4$.   
\item
$\DD$ resolves to a morphism $\widetilde{\DD}$ 
from a blowing-up $\widetilde{\M_0}$,
$$
\begin{array}{ccccc}
\E & \subset & \widetilde{\M_0} &  &  \\
\downarrow & & \downarrow & & \\
\widetilde{\K_0} & \subset & \BL_s(\M_0) & \searrow^{\widetilde{\DD}} & \\
\downarrow & & \downarrow & & \\
\K_0 & \subset & \M_0 & \map{\DD} & \M_\omega  
\end{array}
$$
where $\widetilde{\M_0}$ is obtained by two successive blowing-ups:
first we blow-up the singular points of $\K_0$ and secondly we blow-up 
$\BL_s(\M_0)$ along the smooth proper transform $\widetilde{\K}_0$ of $\K_0$.  
The exceptional divisor $\E$ is mapped by $\widetilde{\DD}$ 
onto the divisor $\T_\omega$.
\end{enumerate}
\end{thm}

\subsection{Restriction of $\DD$ to the extension spaces}

The strategy of the proof is to restrict $\DD$ to the extension spaces
$\pp_0(L)$. We start by defining a map 
$$ \DD_L: \pp_0(L) \lra \M_\omega $$
as follows: consider a point $e \in \pp_0(L)$ \eqref{exte} and denote
by $W_e \subset H^0(C, \omega L^2)$ the corresponding $3$-dimensional
linear subspace of divisors. If we suppose that $e \notin \varphi(C)$,
then the evaluation map $\cO_C \otimes W_e \map{ev} \omega L^2$
is surjective and we define $F_e = \DD_L(e)$ to be the rank $2$
vector bundle such that $\ker(ev) \cong \left( F_e L \right)^*$,i.e.,
we have an exact sequence
\begin{equation} \label{deffe}
 0 \lra \left( F_e L \right)^* \lra \cO_C \otimes W_e \map{ev}
\omega L^2 \lra 0. 
\end{equation}

If there is no ambiguity, we will drop the subscript $e$.
\begin{lem} \label{bundlef}
Suppose that $e \notin \varphi(C)$. Then
\begin{enumerate}
\item The bundle $F_e$ has canonical determinant, is semi-stable and
$F_eL$ is generated by global sections.
\item There exists a nonzero map $L \ra F_e$, hence $[F_e]$ defines a point 
in $\pp_\omega(L)$.
\item We have $\dim H^0(C, E \otimes F_e) \geq 4$, where $E$ is the
stable bundle associated to $e$ \eqref{exte}.
\end{enumerate}
\end{lem}

\begin{proof}
(1) The first assertion is immediately deduced from the exact sequence 
\eqref{deffe}. We take the dual of \eqref{deffe}
\begin{equation} \label{deffe2}
0 \lra \omega^{-1}L^{-2} \lra \cO_C \otimes W^* \lra FL \lra 0.
\end{equation}
Taking global sections leads to the inclusion $W^* \subset H^0(FL)$,
which proves the last assertion. Let us check semi-stability:
suppose that there exists a line subbundle $M$ which 
destabilizes $FL$ (assume $M$ saturated),i.e., $0 \ra M \ra FL \ra
\omega L^2 M^{-1} \ra 0$. Then $\deg M \geq 4$, which implies that
$\deg \omega L^2 M^{-1} \leq 2$. Hence $\dim H^0(\omega L^2 M^{-1})
\leq 1$, so the subspace $H^0(M) \subset H^0(FL)$ has codimension
$\leq 1$, which contradicts that $FL$ is globally generated.

(2) Since $\det F = \omega$, we have $\left( FL \right)^*  = 
FL^{-1} \omega^{-1}$. Taking global sections of the exact
sequence \eqref{deffe} tensored with $\omega$ leads to
$$ 0 \lra H^0(FL^{-1}) \lra H^0(\omega) \otimes W \lra 
H^0(\omega^2L^2) \lra \cdots $$
Now we observe that $\dim H^0(\omega) \otimes W = 9$ and 
$\dim H^0(\omega^2 L^2) = 8$ (Riemann-Roch), which implies
that $\dim H^0(FL^{-1}) \geq 1$.

(3) We tensor the exact sequence \eqref{exte} defined by $e$ with $F$
and take global sections
$$ 0 \lra H^0(FL^{-1}) \lra H^0(E \otimes F) \lra H^0(FL) 
\map{\cup e} H^1(FL^{-1}) \lra \cdots $$
The coboundary map is the cup-product with the extension class
$e \in H^1(L^{-2})$ and, since $\det F = \omega$, the coboundary
map $\cup e$ is skew-symmetric (by Serre-duality $H^1(FL^{-1})= 
H^0(FL)^*$). Hence the linear map $\epsilon \mapsto \cup \epsilon$
factorizes as follows
\begin{equation} \label{cupp}
H^0(\omega L^2)^* \lra \Lambda^2 H^0(FL)^* \subset 
\Hom(H^0(FL), H^1(FL)),
\end{equation} 
and its dual map $\Lambda^2 H^0(FL) \map{\mu} H^0(\omega L^2)$
coincides with exterior product of global sections (see e.g. \cite{la1}).
On the other hand, it is easy to check that the image under $\mu$
of the subspace $\Lambda^2 W^* \subset \Lambda^2 H^0(FL)$ equals
$W \subset H^0(\omega L^2)$ and $\mu$ restricts to the canonical
isomorphism $\Lambda^2 W^* = W$. Therefore the linear map  
$\cup e$ is zero on $W^* \subset H^0(FL)$, from which we deduce 
that $\dim H^0(E \otimes F) = \dim H^0(FL^{-1}) + \dim \ker(\cup e)
\geq 4$.
\end{proof}

It follows that the map $\DD_L$ factorizes  
\begin{equation} \label{polarext}
\DD_L : \pp_0(L) \lra \pp_\omega(L) \subset \M_\omega.
\end{equation}
Moreover, by Lemma \ref{bundlef}(3) and Lemma \ref{tansing} 
the point $\DD_L(e)$ corresponds to the
embedded tangent space at $e \in \pp_0(L)$, hence $\DD_L$ is the
restriction of $\DD$ to $\pp_0(L)$. In particular, $\DD_L$ is given
by a linear system of cubics through $\varphi(C)$.

We recall that the restriction of the trisecant scroll $\T_0$ to 
$\pp_0(L)$ is the
surface, denoted by $\T_0(L)$, ruled out by the trisecants to
$\varphi(C) \subset \pp_0(L)$.

\begin{lem} \label{trisecant}
Given a point $e \in \pp_0(L)$ such that $e \notin \varphi(C)$. The bundle
$F_e$ is stable if and only if $e \notin \T_0$. Moreover 
\begin{itemize}
\item if $\dim H^0(L^2) = 0$, then the 
trisecant $\overline{pqr}$ to $\varphi(C)$ is contracted to the 
semi-stable point $[L(u) \oplus \omega L^{-1}(-u)] =
\varphi(u) \in \pp_\omega(L)$ for some point $u \in C$, which 
satisfies $p+q+r \in |L^2(u)|$.
\item if $\dim H^0(L^2) >0$, then $\omega L^{-2} = \cO_C(u+v)$ for some
points $u,v \in C$ and any trisecant $\overline{pqr}$ is contracted to
the semi-stable point $[L(u) \oplus L(v)]$.
\end{itemize}
\end{lem}

\begin{proof}
The bundle $F$ fits into an exact sequence $0 \ra L \ra F \ra \omega L^{-1}
\ra 0$. Suppose that $F$ has a line subbundle $M$ of degree $2$ and
consider the composite map $\alpha: M \ra F \ra \omega L^{-1}$.

\bigskip
\noindent
First we consider the case $\alpha = 0$: then $M = L(u) \hookrightarrow
F$ for some $u \in C$, or equivalently $\dim H^0(FL^{-1}(-u))>0$. 
We tensor \eqref{deffe} with $\omega(-u)$ and take global sections
$$ 0 \lra H^0(FL^{-1}(-u)) \lra H^0(\omega(-u)) \otimes W \map{m}
H^0(\omega^2L^2(-u)) \lra \cdots $$
The second map $m$ is the multiplication map of global 
sections. Since $W \subset H^0(\omega L^2)$, let us 
consider for a moment the extended multiplication map
$\tilde{m}: H^0(\omega(-u)) \otimes H^0(\omega L^2) \lra
H^0(\omega^2L^2(-u))$. By the base-point-free-pencil-trick applied
to the pencil $|\omega(-u)|$, we have $\ker \tilde{m}  = H^0(L^2(u))$
and a tensor in $\ker \tilde{m}$ is of the form $s \otimes t\alpha -
t \otimes s \alpha$, with $\{s,t\}$ a basis of $H^0(\omega(-u))$
and $\alpha \in H^0(L^2(u))$. We denote by $p+q+r$ the zero divisor
of $\alpha$. Then we see that $\ker m \not= \{ 0 \}$ if and
only if $W$ contains the linear space spanned by $t\alpha$ and
$s\alpha$. Dually, this means that $e \in \overline{pqr}$, the trisecant
through the points $p,q,r$. Conversely, any $e \in \overline{pqr}$ 
is mapped by $\DD_L$ to $[L(u) \oplus \omega L^{-1}(-u)]$.

\bigskip
\noindent
Secondly we consider the case $\alpha \not= 0$: then $M = \omega L^{-1}(-u)
\hookrightarrow F$ for some $u \in C$, or equivalently $\dim
H^0(F \omega^{-1}L(u)) >0$. As in the first case we take global
sections of \eqref{deffe} tensored with $L^2(u)$ and we obtain that
$H^0(F\omega^{-1}L(u))$ is the kernel of the multiplication
map $H^0(L^2(u)) \otimes W \map{m} H^0(\omega L^4(u))$. Then
$\ker \tilde{m} \not= \{0 \}$ implies that $\dim H^0(L^2(u)) = 2$. 
Hence $L^2(u) = \omega(-v)$ for some point $v \in C$,i.e., $\omega L^{-2}
= \cO_C(u+v)$, which implies that $\dim H^0(\omega L^{-2}) = 
\dim H^0(L^2) > 0$. Furthermore the multiplication map
becomes $H^0(\omega(-v)) \otimes W \map{m} H^0(\omega^2 L^2(-v))$.
We can now conclude exactly as in the first case, with the additional
observation that any trisecant $\overline{pqr}$ is contracted to the
point $[L(v) \oplus \omega L^{-1}(-v)] = [L(v) \oplus L(u)]$. 
\end{proof}

Now we are going to construct along the same lines an inverse map
to $\DD_L$ \eqref{polarext}
$$ \DD'_L : \pp_\omega(L) \lra \pp_0(L). $$
Given an extension class $f \in \PP_\omega(L)$ such that $f \notin
\varphi(C)$, we denote by $W_f \subset H^0(C, \omega^2 L^{-2})$ the
corresponding $3$-dimensional linear space of divisors and we define
$E_f = \DD'_L(f)$ to be the rank $2$ vector bundle which fits
in the exact sequence
$$ 0 \lra E_f \omega^{-1}L \lra W_f \otimes \cO_C \map{ev} 
\omega^2 L^{-2} \lra 0. $$
Exactly as in Lemma \ref{bundlef} we show that $E_f$ has the
following properties.

\begin{lem} \label{bundlefdual}
Suppose that $f \notin \varphi(C)$. Then
\begin{enumerate}
\item The bundle $E_f$ has trivial determinant, is semi-stable and
$E_f\omega L^{-1}$ is generated by global sections.
\item There exists a nonzero map $L^{-1} \ra E_f$, hence $[E_f]$ 
defines a point in $\pp_0(L)$.
\item We have $\dim H^0(C, E_f \otimes F) \geq 4$, where $F$ is the
stable bundle associated to $f$ \eqref{extf}.
\end{enumerate}
\end{lem}

Similarly the analogue of Lemma \ref{trisecant} holds for the bundle $E_f$.

\begin{lem}
The map $\DD'_L$ is the birational inverse of $\DD_L$,i.e.,
$$ \DD'_L \circ \DD_L = Id_{\pp_0(L)} \qquad \text{and} 
\qquad \DD_L \circ \DD'_L = Id_{\pp_\omega(L)}. $$
\end{lem}

\begin{proof}
Start with $e \in \pp_0(L)$ with $e \notin \T_0(L)$. Then 
(Lemma \ref{trisecant}) $\DD_L(e) = F_e$ is stable and (Lemma
\ref{bundlef}(3)) $\dim H^0(C, E \otimes F_e) \geq 4$. Now the
stable bundle $F_e$ determines an extension class $f \in \pp_\omega(L)$
with $f \notin \varphi(C)$. Let us denote $E_f = \DD'_L(f)$. We know
(Lemma \ref{bundlefdual}(3)) that $\dim H^0(C, E_f \otimes F_e) \geq 4$
and since $F$ is stable we deduce from Lemma \ref{tansingdual}
that the embedded tangent space $\TT_F \M_\omega$ corresponds
to $[E]$ and $[E_f]$. Hence $[E]=[E_f]$ and since $E$ is stable, we
have $E = E_f$.
\end{proof}

We deduce that $\DD_L$ restricts to an isomorphism
$\pp_0(L) \setminus \T_0(L) \map{\sim} \pp_\omega(L) \setminus \T_\omega(L)$.
Since $\M_0$ is covered by the spaces $\pp_0(L)$ and since $\DD$ restricts
to $\DD_L$ on $\pp_0(L)$, we easily deduce from the preceeding lemmas
parts (1),(2) and (3) of Theorem \ref{maintheorem}.

\subsection{Blowing-up}

Even if part (4) of Theorem \ref{maintheorem} is a straight-forward
consequence of the results obtained in \cite{la2}, we give the 
complete proof for the convenience of the reader.
First we consider the blowing-up $\BL_s(\pp^7)$ of $\pp^7 = |2\Theta|$
along the $64$ singular points of $\K_0$. Because of the
invariance of $\K_0$ and $\M_0$ under the Heisenberg group, it is
enough to consider the blowing-up at the ``origin'' $O := [\cO \oplus
\cO]$. We denote by 
$\widetilde{\K_0}$ (resp. $\BL_s(\M_0)$) the proper transform
of $\K_0$ (resp. $\M_0$), and by $\pp(T_O \pp^7) \subset \BL_s(\pp^7)$
the exceptional divisor (over $O$).

\bigskip
\noindent
By \cite{la2} Remark 5 the Zariski tangent spaces $T_O\K_0$ and
$T_O \M_0$ at the origin $O$ to $\K_0$ and $\M_0$ satisfy the
relations
$$ \sym^2 H^0(\omega)^* \cong T_O \K_0 \subset T_O \M_0 =
T_O \pp^7 \qquad \text{and} \qquad T_O\M_0 / T_O \K_0 \cong 
\Lambda^3 H^0(\omega)^*.$$
Moreover with the notation of section 2.2 the equation of the
hyperplane $T_O \K_0 \subset T_O \M_0$ is $T=0$ and $T_{ij}$ are
coordinates on $\sym^2 H^0(\omega)^*$. We deduce from the
local equation of $\M_0$ at the origin $O$ (section 2.2(ii)) that
$\widetilde{\K}_0 \cap \pp \sym^2 H^0(\omega)^*$ is the Veronese
surface $S := \mathrm{Ver} H^0(\omega)^*$ and that $\widetilde{\K}_0$
is smooth. Moreover the linear system spanned by the proper
transforms of the cubics $C_i$ is given by the six quadrics
$Q_{ij} := \frac{\partial}{\partial T_{ij}}\left( \det [T_{ij}] \right)$
vanishing on $S$.

\bigskip
\noindent

Given a smooth point $x= [M \oplus M^{-1}] \in \K_0$ with $M^2 \not= \cO$,
the Zariski tangent spaces $T_x\K_0$ and $T_x\M_0$ satisfy the relations
$$ H^0(\omega)^* \cong T_x \K_0 \subset T_x \M_0 = T_x \pp^7 
\qquad \text{and} \qquad
T_x \M_0 / T_x \K_0 \cong H^0(\omega M^2)^* \otimes H^0(\omega M^{-2})^*.$$
The tangent space $T_x \K_0 \subset T_x \M_0$ is cut out by the
four equations $T_{ij} = 0$, where the $T_{ij}$ are natural coordinates
on $H^0(\omega M^2)^* \otimes H^0(\omega M^{-2})^*$. Let $\widetilde{\E}$
be the exceptional divisor of the blowing-up of $\BL_s(\pp^7)$ along
the smooth variety $\widetilde{\K_0}$ and let $\E$ be its restriction
to the proper transform $\widetilde{\M_0}$. We denote by $\widetilde{\E}_x$
and $\E_x$ the fibres of $\widetilde{\E}$ and $\E$ over a point $x \in
\K_0$. Then for a smooth point $x$, it follows from the local equation
at $x$ (section 2.2 (ii)) that  $\E_x$ is the Segre embedding
$\pp^1 \times \pp^1 = |\omega M^2|^* \times |\omega M^{-2}|^*
\hookrightarrow \pp H^0(\omega M^2)^* \otimes H^0(\omega M^{-2})^* =
\widetilde{\E}_x$ and the linear system spanned by the proper
transforms of the cubics $C_i$ is given by the four linear forms
$T_{ij}$.

\bigskip
\noindent

At a singular point (we take $x = O$), it follows from the
preceeding discussion that $\E_O$ is the exceptional divisor 
of the blowing-up of $\pp \sym^2 H^0(\omega)^*$ along the
Veronese surface $S$,i.e., the projectivized normal bundle
over $S$. It is a well-known fact (duality of conics) that 
the rational map given by the quadrics $Q_{ij}$ resolves
by blowing-up $S$.

\bigskip
\noindent

It remains to show that $\widetilde{\DD}$ maps $\E$ onto the
trisecant scroll $\T_\omega$. Since $\E$ is irreducible, it will
be enough to check this on an open subset of $\E$. We consider
again the extension spaces $\pp_0(L) \subset \M_0$. For simplicity
we choose $L$ such that
\begin{itemize}
\item[(1)] $\pp_0(L)$ does not contain a singular point of $\K_0$,
\item[(2)] the morphism $\varphi: C \lra \pp_0(L)$ is an embedding, or
equivalently $\dim H^0(L^2) = 0$.
\end{itemize} 
Let $\widetilde{\pp_0(L)}$ be the blowing-up of $\pp_0(L)$ along the
curve $C$, with exceptional divisor $\E_L$. Because of assumptions (1) and
(2), we have an embedding $\widetilde{\pp_0(L)} \hookrightarrow 
\widetilde{\M_0}$, $\E$ restricts to $\E_L$, and $\E_L$ is the projectivized
normal bundle $N$ of the embedded curve $C \subset \pp_0(L)$. We have the
following commutative diagram
$$
\begin{array}{ccccc}
\pp(N) = \E_L & \subset & \widetilde{\pp_0(L)} & & \\
 & & & & \\
\downarrow^{\pi} & & \downarrow & \searrow^{\widetilde{\DD}_L} & \\
 & & & & \\
C & \subset & \pp_0(L) & \map{\DD_L} & \pp_\omega(L)
\end{array}
$$
In order to study the image $\widetilde{\DD}_L(\E_L)$ we introduce,
for a point $u \in C$, the rank $2$ bundle $E_u$ which
is defined by the exact sequence
$$ 0 \lra E^*_u \lra \cO_C \otimes H^0(\omega L^2(-u)) \map{ev}
\omega L^2(-u) \lra 0.$$
Note that $H^0(\omega L^2(-u))$ corresponds to the
hyperplane defined by $u \in C \subset \pp_0(L)$. Then exactly as
in Lemma \ref{bundlef}(1) we show that $\det E_u = \omega L^2(-u)$,
$E_u$ is stable and globally generated with $H^0(E_u) \cong
H^0(\omega L^2(-u))^*$. We introduce the Hecke line $\H_u$ defined
as the set of bundles which are (negative) elementary transformations
of $E_uL^{-1}(u)$ at the point $u$,i.e., the set of bundles which
fit into the exact sequence 
\begin{equation} 
0 \lra F \lra E_uL^{-1}(u) \lra \CC_u \lra 0.
\end{equation} 
Since $E_u$ is stable, any $F$ is semi-stable (and $\det F = \omega$)
and we have a linear map (see \cite{b2}) $\pp^1 \cong \H_u
\ra \M_\omega$.

\begin{lem} \label{heckeline}
Given a point $u \in C$, the fibre $\pp(N_u) = \E_{L,u}$ is mapped
by $\widetilde{\DD}_L$ to the Hecke line $\H_u \subset \pp_\omega(L)$.
Moreover $\H_u$ coincides with the trisecant line $\overline{pqr}$
to $C \subset \pp_\omega(L)$, with $p+q+r \in |\omega L^{-2}(u)|$.
\end{lem}

\begin{proof}
The Zariski tangent space $T_u \pp_0(L)$ at the point $u$ is identified
with $H^0(\omega L^2(-u))^* \cong H^0(E_u)$. Under this identification
the tangent space $T_u C$ corresponds to the subspace $H^0(E_u(-u))$. 
Hence we obtain a canonical isomorphism of $\pp(N_u)$ with the
projectivized fibre over the point $u$ of the bundle $E_u$,i.e.,
the Hecke line $\H_u$. It is straight-forward to check that
$\widetilde{\DD}_L$ restricts to the isomorphism $\pp(N_u) \cong
\H_u$. To show the last assertion, it is enough to observe
that the Hecke line $\H_u$ intersects the curve 
$C \subset \pp_\omega(L)$ at a point $p$ if and only if 
$\dim H^0(E_u L^{-1}(u-p))>0$ and to continue as in the 
proof of Lemma \ref{trisecant}.
\end{proof}

Since the union of those $\E_L$ such that $L$ satisfies assumptions (1) and (2)
form an open subset of $\E$, we conclude that 
$\widetilde{\DD}(\E) = \T_\omega$.
This completes the proof of Theorem \ref{maintheorem}.

\subsection{Some remarks} 

{\bf{(1)}}
The divisor $\T_\omega \in |\LL^8|$. This is seen as follows: it suffices
to restrict $\T_\omega$ to a general $\pp_\omega(L) \subset \M_\omega$
and to compute the degree of the trisecant scroll $\T_\omega(L) \subset
\pp_\omega(L)$. By Lemma \ref{heckeline} $\T_\omega(L)$ is the image
of $\E_L = \pp(N)$ under the morphism $\widetilde{\DD}_L$. The
hyperplane bundle over $\pp_\omega(L)$ pulls-back under $\widetilde{\DD}_L$
to $\cO_{\pp}(1) \otimes \pi^*\left( \omega^3 L^6 \right)$ over the
ruled surface $\pp(N)$. Since $\widetilde{\DD}_{L | {\E}_L}$ is 
birational, we
obtain that $\deg \T_\omega(L) = \deg  \pi_* \cO_{\pp}(1) \otimes
\omega^3 L^6 = \deg N^* \omega^3 L^6 = 8$.

\bigskip
\noindent
{\bf{(2)}}
Using the same methods as before, one can show a refinement
of Theorem \ref{maintheorem}(3). Consider $E$ stable with $E
\in \M_0$ and $F$ semi-stable with $[F] \in \M_\omega$. 
\begin{itemize}
\item
The only pairs $(E,F)$ for which $\dim H^0(C, E\otimes F) =  6$ are 
the $64$ exceptional pairs $E = A_\kappa$ and $F = \kappa \oplus \kappa$
for a theta-characteristic $\kappa$ \eqref{distbun}. We note that
$\DD(A_\kappa) = [\kappa \oplus \kappa]$.

\item Suppose $\DD(E) = [M \oplus \omega M^{-1}]$ for some $M$ and
$E \not= A_\kappa$,i.e.,$M^2 \not= \omega$. 
Then there are exactly three semi-stable 
bundles $F$ such that $\DD(E) = [F]$ and 
$\dim H^0(C, E \otimes F) = 4$, namely
\begin{itemize}
\item[(1)] the decomposable bundle $F = M \oplus \omega M^{-1}$ (note that
$\dim H^0(EM) = 2$).
\item[(2)] two indecomposable bundles with extension classes in
$\ext^1(M,\omega M^{-1}) = H^0(M^2)^*$ and 
$\ext^1(\omega M^{-1},M) = H^0(\omega^2 M^{-2})^*$ defined by
the images of the exterior product maps
$$ \Lambda^2 H^0(EM) \lra H^0(M^2) \qquad \text{and}
\qquad \Lambda^2 H^0(E\omega M^{-1}) \lra H^0(\omega^2 M^{-2}). $$
\end{itemize} 
\end{itemize}

\bigskip
\noindent
{\bf{(3)}}
As a corollary of Lemma \ref{heckeline} we obtain that
the morphism $\widetilde{\DD}$ maps the exceptional
divisor $\widetilde{\E}$ onto the dual hypersurface $\K^*_0$ of the Kummer
variety $\K_0$ (more precisely $\widetilde{\DD}$ maps
$\widetilde{\E}_x = \pp^3$ isomorphically to the subsystem of
divisors singular at $x \in \K_0^{sm}$) 
and that the hypersurface $\widetilde{\DD}(\widetilde{\E})
= \K^*_0$ intersects (set-theoretically) $\M_\omega$ along the
trisecant scroll $\T_\omega$. It is worthwhile to figure out 
the relationship with other distinguished hypersurfaces in $|2\Theta|$,
e.g. the octic $G_8$ defined by the equation $\DD^{-1}(F_4) = F_4 \cdot
G_8$ and the Hessian $H_{16}$ of Coble's quartic $F_4$.

\section{Applications}

\subsection{The $8$ maximal line subbundles of $E \in \M_0$}

In this section we recall the results of \cite{ln} (see also
\cite{opp}, \cite{op}) on line subbundles of stable bundles
$E \in \M_0$ and $F \in \M_\omega$. We introduce the
closed subsets $\MM_0(E)$ and $\MM_\omega(F)$ of $\pic^1(C)$
parametrizing  line subbundles of maximal degree of $E$ and
$F$,
$$ \MM_0(E) := \{ L \in \pic^1(C) \ | \ L^{-1} \hookrightarrow
E \} \qquad \text{and} \qquad 
\MM_\omega(F) := \{ L \in \pic^1(C) \ | \ L \hookrightarrow
F \}.$$
The next lemma follows from \cite{ln} section 5 and Nagata's theorem. For
simplicity we assume that $C$ is not bi-elliptic.

\begin{lem}
The subsets $\MM_0(E)$ and $\MM_\omega(F)$ are non-empty and
$0$-dimensional, unless $E$ and $F$ are exceptional  
(see \eqref{distbun}). In these cases we have
$$\MM_0(A_\kappa) = \{ \kappa(-p) \ | \ p \in C \} \cong C 
\qquad \text{and} \qquad
\MM_\omega(A_\alpha) = \{ \alpha(p) \ | \ p \in C \} \cong C. $$
\end{lem}

Note that $A_\kappa \in \T_0$ and $A_\alpha \in \T_\omega$ (see \cite{opp}
Theorem 5.3) and that in the bi-elliptic case, we additionally
have a $JC[2]$-orbit in $\M_0$ (resp. $\M_\omega$) of bundles
$E$ (resp. $F$) with $1$-dimensional $\MM_0(E)$ (resp.
$\MM_0(F)$).

\bigskip
\noindent
Since $\MM_0(E)$ is non-empty, any stable $E \in \M_0$ lies in at least
one extension space $\pp_0(L)$ for some $L \in \pic^1(C)$
with extension class $e \notin \varphi(C)$. Now \cite{ln}
Proposition 2.4 says that there exists a bijection between the sets of
\begin{itemize}
\item[(1)] effective divisors $p+q$ on $C$ such that $e$ lies on
the secant line $\overline{pq}$
\item[(2)] line bundles $M \in \pic^1(C)$ such that $M^{-1} \hookrightarrow
E$ and $M \not= L$.
\end{itemize} 
The two data are related by the equation
\begin{equation} \label{rellsb}
L \otimes M = \cO_C(p+q).
\end{equation}
Let us count secant lines to $\varphi(C)$ through a {\em general} point 
$e \in \pp_0(L)$: composing $\varphi$ with the projection away
from $e$ maps $C$ birationally to a plane nodal sextic $S$. By the
genus formula, we obtain that the number of 
nodes of $S$ (= number
of secants) equals $7$. Hence, for $E$ general, the cardinal
$|\MM_0(E)|$
of the finite set $\MM_0(E)$ is $8$. We write
$$\MM_0(E) = \{ L_1, \cdots, L_8 \}.$$

From now on we assume that $E$ is sufficiently general in order to
have $|\MM_0(E)| = 8$. Since $E \in \pp_0(L_i)$ for $1 \leq i \leq 8$,
we deduce from relation \eqref{rellsb} that
\begin{equation} \label{gbt}
1 \leq i < j \leq 8, \qquad L_i \otimes L_j = \cO_C(D_{ij}),
\end{equation}
where $D_{ij}$ is an effective degree two divisor on $C$.

\begin{lem}
The $8$ line bundles $L_i$ satisfy the relation $\bigotimes_{i=1}^8 L_i
= \omega^2$.
\end{lem}

\begin{proof}
We represent $E$ as a point $e \in \pp_0(L_8)$ and assume
that the plane sextic curve $S \subset \pp^2$ obtained by
projection with center $e$ has $7$ nodes as singularities. 
It will be enough to prove the equality for such a bundle $E$.
Then $C \map{\pi} S$ is the normalization of $S$ and, by the
adjunction formula, we have $\omega = \pi^* \cO_{S}(3) \otimes
\cO_C(-\Delta)$, where $\Delta$ is the divisor lying over the $7$ nodes of
$S$,i.e., $\Delta = \sum_{i=1}^{7} D_{i8}$. Hence
$$\omega = \omega^3 L_8^6 \left(- \sum_{i=1}^{7} D_{i8}\right) 
= \omega^3 L_8^{-1} \otimes \bigotimes_{i=1}^{7} \left( L_8 (-D_{i8})
\right) = \omega^3 \otimes \bigotimes_{i=1}^{8} L_i^{-1},$$
where we used relations \eqref{gbt}.
\end{proof}

\begin{rem} \label{reconst}
Conversely, suppose we are given $8$ line bundles $L_i$ which
satisfy the $28$ relations \eqref{gbt}. Then there exists
a unique stable bundle $E \in \M_0$ such that $\MM_0(E) =
\{ L_1, \ldots, L_8 \}$. This is seen as follows: take e.g.
$L_8$ and consider any two secant lines $\overline{D}_{i8}$ and
$\overline{D}_{j8}$ ($i<j<8$) in $\pp_0(L_8)$. Then relations
\eqref{gbt} imply that these two lines intersect in a point $e$.
It is straight-forward to check that the bundle $E$ associated 
to $e$ does not depend on the choices we made.
\end{rem}

\subsection{Nets of quadrics}

We consider 
$E \in \M_0$ and we assume that $E \notin \T_0$ and $|\MM_0(E)| = 8$.
Then $F = \DD(E)$ is stable and $\dim H^0(C, E \otimes F) = 4$.
We recall that the rank $4$ vector bundle $E \otimes F$ is equipped
with a non-degenerate quadratic form 
$$ \det : E \otimes F = \Hom (E,F) \lra \omega. $$
Taking global sections on both sides endows the projective
space $\pp^3 := \pp H^0(C, \Hom (E,F))$ with a net $\Pi = |\omega|^*$ of 
quadrics. We denote by $Q_x \subset \pp^3$ the quadric
associated to $x \in \Pi$ and, identifying $C$ with its
canonical embedding $C \subset |\omega|^* = \Pi$, we see that (the cone
over) the quadric $Q_p$ for $p \in C$ corresponds to the sections 
\begin{equation} \label{quadric}
 Q_p := \{ \phi \in H^0(C, \Hom(E,F)) \ | \ E_p \map{\phi_p} F_p \ 
\text{not surjective} \}, 
\end{equation}
where $E_p, F_p$ denote the fibres of $E,F$ over $p \in C$. It
follows from Lemma \ref{bundlef}(2) that $\MM_0(E) = \MM_\omega(F)$,or
equivalently any line bundle $L_i \in \MM_0(E)$ fits into a 
sequence of maps $$x_i: E \lra  L_i \lra F.$$
We denote by $x_i \in \pp^3$ the composite map (defined up to a scalar).

\begin{lem}
The base locus of the net of quadrics $\Pi$ consists of the $8$
distinct points $x_i  \in \pp^3$.
\end{lem}

\begin{proof}
A base point $x$ corresponds to a vector bundle map $x: E \ra F$ such that
$\rk x \leq 1$ (since $x \in Q_p$, $\forall p$). Hence there exists
a line bundle $L$ such that $E \ra L \ra F$ and since $E$ and $F$ are
stable, of slope $0$ and $2$, we obtain that $\deg L = 1$ and $L \in 
\MM_0(E) = \MM_\omega(F)$.
\end{proof}

The set of base points $\overline{x} = \{ x_1,\ldots,x_8 \}$ of a net
of quadrics in $\pp^3$ is {\em self-associated} (for the definition
of (self-)association of point sets we refer to \cite{do} chapter 3)
and called a {\em Cayley octad}. We recall (\cite{do} chapter 3 example
6) that ordered Cayley octads $\overline{x} = \{ x_1,\ldots,x_8 \}$
are in 1-to-1 correspondence with ordered point sets $\overline{y} =
\{ y_1,\ldots,y_7 \}$ in $\pp^2$ (note that we consider here 
general ordered point sets up to  projective equivalence).
The correspondence goes as follows: starting from $\overline{x}$
we consider the projection with center $x_8$, $\pp^3 \map{pr_{x_8}} \pp^2$,
and define $\overline{y}$ to be the projection of the remaining seven points.
Conversely, given $\overline{y}$ in $\pp^2$, we obtain by association
$7$ points $x_1,\ldots,x_7$ in $\pp^3$. The missing $8$-th point $x_8$
of $\overline{x}$ is the additional base point of the net of quadrics
through the $7$ points $x_1,\ldots,x_7$.

\bigskip
\noindent
Consider a general $E \in \M_0$ and choose a line 
subbundle $L_8 \in \MM_0(E)$. We denote by $x_8$ the 
corresponding base point of the net $\Pi$. We consider the
two (different) projections onto $\pp^2$
\begin{itemize}
\item[(1)] projection with center $x_8$ of $\pp^3 = \pp H^0(C, \Hom (E,F))
\map{pr_{x_8}} \pp^2$. Let $\overline{y} = \{y_1, \ldots, y_7 \} \subset
\pp^2$ be the projection of the $7$ base points $x_1, \ldots, x_7$.
\item[(2)] projection with center $e$ of $\pp_0(L_8) \map{pr_e} \pp^2$.
Let $\overline{z} = \{z_1, \ldots , z_7 \} \subset \pp^2$ be the images
of the $7$ secant lines to $\varphi(C)$ through $e$. Note that
$z_1,\ldots, z_7$ are the $7$ nodes of the plane sextic $S$. 
\end{itemize}

\begin{lem} \label{ident}
The two target $\pp^2$'s of the projections (1) and (2) are canonically
isomorphic (to $\pp W_e^*$) and the two point sets $\overline{y}$ and
$\overline{z}$ coincide.
\end{lem}

\begin{proof}
First we recall from the proof of Lemma \ref{bundlef} that we have an
exact sequence 
$$ 0 \lra H^0(FL^{-1}_8) \map{i} H^0(E \otimes F) \map{\pi} 
H^0(FL_8) \lra 0, $$
and that $H^0(FL_8) \cong W_e^*$ and $\dim H^0(FL_8^{-1}) =1$.
Moreover it is easily seen that $\pp(\im i) = x_8 \in \pp^3$, hence
the projectivized map $\pi$ identifies with $pr_{x_8}$. The images
$pr_{x_8}(x_i)$ for $1 \leq i \leq 7$ are given by the sections
$s_i \in H^0(FL_8)$ vanishing at the divisor $D_{i8}$
(since $L_iL_8 = \cO_C(D_{i8}) \hookrightarrow FL_8$). It remains
to check that the section $s_i \in H^0(FL_8) \cong W_e^*$
correponds to the $2$-dimensional subspace $H^0(\omega L^2(-D_{i8}))
\subset W_e \subset H^0(\omega L^2)$, which is standard. 
\end{proof}

We introduce the non-empty open subset $\M_0^{reg} \subset \M_0$ of
stable bundles $E$ which satisfy $E \notin \T_0$, $|\MM_0(E)| = 8$ and
for any $L \in \MM_0(E)$ the point set $\overline{z} \subset \pp^2$
is such that no three points in $\overline{z}$ are collinear.

\subsection{The Hessian construction}

It is classical (see e.g. \cite{do} chapter 9) 
to associate to a net of quadrics $\Pi$ on $\pp^3$ its Hessian curve 
parametrizing singular quadrics,i.e.,
$$ \Hess(E) := \{ x \in \Pi = |\omega|^*  \ | \ Q_x  \ \text{singular} \}. $$
Note that $C$ and $\Hess(E)$ lie in the same projective plane.

\begin{lem}
We suppose that $E \in \M_0^{reg}$. Then the curve $\Hess(E)$ is a 
smooth plane quartic.
\end{lem}

\begin{proof}
It follows from \cite{do} chapter 9 Lemma 5 that $\Hess(E)$ is smooth
if and only if every $4$ points of $\overline{x} = \{ x_1,\ldots,x_8 \}$
span $\pp^3$. Projecting away from one of the $x_i$'s  and using Lemma
\ref{ident} we see that this condition holds for $E \in \M_0^{reg}$.
\end{proof}

First we determine for which bundles $E \in \M_0^{reg}$ the Hessian
curve $\Hess(E)$ equals the base curve $C$. We need to recall some
facts about nets of quadrics and Cayley octads \cite{do}. 
The net $\Pi$ determines
an even theta-characteristic $\theta$ over the smooth curve 
$\Hess(E)$, such that the Steinerian embedding
$$ \Hess(E) \map{\St} \pp^3 = |\omega \theta|^*, \qquad
x \lms \Sing (Q_x), $$
is given by the complete linear system $|\omega \theta|$. The image
$\St(E)$ is called the {\em Steinerian} curve. Given two distinct
base points $x_i,x_j \in \pp^3$ of the net $\Pi$, the pencil
$\Lambda_{ij}$ of quadrics of the net $\Pi$ which contain
the line $\overline{x_i x_j}$ is a bitangent to the curve
$\Hess(E)$. In this way we obtain all the $28 = \binom{8}{2}$ bitangents
to $\Hess(E)$. Let $u,v$ be the two intersection points of the
bitangent $\Lambda_{ij}$ with $\Hess(E)$. Then the secant line 
to the Steinerian curve $\St(E)$ determined by $\St(u)$ and $\St(v)$ 
coincides with $\overline{x_i x_j}$.

\bigskip
\noindent

Conversely, given a smooth plane quartic $X \subset \pp^2$ with an
even theta characteristic $\theta$, taking the symmetric resolution 
over $\pp^2$ of the
sheaf $\theta$ supported at the curve $X$  
gives a net of quadrics $\Pi$ whose Hessian curve equals $X$.
Thus the correspondence between nets of quadrics $\Pi$ and 
the data $(X,\theta)$ is 1-to-1.

\bigskip
\noindent

This correspondence allows us to construct some more distinguished
bundles in $\M_0$. We consider a triple $(\theta,L,x)$ consisting of an even 
theta-characteristic $\theta$ over $C$, a square-root $L \in \pic^1(C)$,i.e.,
$L^2 = \theta$, and a base point $x$ of the net of
quadrics $\Pi$ associated to $(C, \theta)$. We denote by
\begin{equation} \label{Aronholdbun}
A(\theta,L,x) \in \M_0 
\end{equation}
the stable bundle defined by the point $x \in \pp_0(L) = |\omega \theta|^*$.
Since $C$ is smooth, we have $A(\theta,L,x) \in \M_0^{reg}$. These bundles
will be called {\em Aronhold} bundles (see Remark \ref{remaronhold}). We leave it
to the reader to deduce the following characterization: 
$E$ is an Aronhold bundle if and only if the $28$ line bundles $L_i L_j$
($1 \leq i<j \leq 8$) are the 
odd theta-characteristics, with $L_i \in \MM_0(E)$.

\begin{prop} \label{evtan}
Given a bundle $E \in \M_0^{reg}$. Then
\begin{enumerate}
\item We have $\Hess(E) = C$ if and only if $E$ is an Aronhold bundle.
\item Assuming $\Hess(E) \not= C$, the curves $C$ and $\Hess(E)$ 
are everywhere tangent. More 
precisely, the scheme-theoretical intersection $C \cap \Hess(E)$
is non-reduced of the form $2\Delta(E)$, with 
$\Delta(E) \in |\omega^2|$.
\end{enumerate}
\end{prop}

\begin{proof}
We deduce from \eqref{quadric} that the intersection $C \cap \Hess(E)$ 
corresponds (set-theoretically) to the sets of points where the 
evaluation map of global sections 
\begin{equation} \label{eval}
\cO_C \otimes H^0(C,\Hom(E,F)) \map{ev} \Hom(E,F)
\end{equation}
is not surjective.

Let us suppose that $C = \Hess(E)$. Then $ev$ is not generically surjective
($\rk ev \leq 3$). We choose a line subbundle $L_8 \in \MM_0(E)$ and
we consider (as in Lemma \ref{ident}) the exact sequence 
$$
\begin{array}{ccccccc}
0 \lra & H^0(FL_8^{-1}) & \lra & H^0(\Hom(E,F)) & \lra & H^0(FL_8) &
\lra 0 \\
 & & & & & & \\
 & \downarrow^{\cong} &   & \downarrow^{ev} &  & \downarrow^{ev'} & \\
 & & & & & & \\
0 \lra & \cO_C & \lra & \Hom(E,F) & \lra & \E & \lra 0
\end{array}
$$
where the vertical arrows are evaluation maps. Note that $\cO_C
\hookrightarrow FL_8^{-1} \hookrightarrow \Hom(E,F)$ corresponds to
the section of $H^0(FL_8^{-1})$. We denote by $\E$ the rank $3$
quotient. Then $ev' : H^0(FL_8) \lra \E$ is not generically 
surjective either. But $\E$ has a quotient $E \ra FL_8$ with 
kernel $\omega L_8^{-2}$. Now since $H^0(FL_8) \map{ev} FL_8$ is
surjective, we obtain a direct sum decomposition $\E = \omega L_8^{-2}
\oplus FL_8$. Furthermore since $E \otimes F$ is poly-stable
(semi-stable and orthogonal) and of slope $2$, we obtain
that $\omega L_8^{-2}$ is an orthogonal direct summand. Hence
$\omega L_8^{-2} = \theta$ for some theta-characteristic $\theta$. 
Now we can do this reasoning for any line bundle $L_i \in \MM_0(E)$,
establishing that all $\omega L_i^{-2}$ are theta-characteristics
contained in $\Hom(E,F)$.  Projecting to $FL_8$ shows that $L_i^2 =  
L_8^2 = \theta$ for all $i$ and therefore the $28$ line bundles
$L_iL_j$ are the odd theta-characteristics. It follows that $E$ is an
Aronhold bundle.

\bigskip

Assuming $C \not= \Hess(E)$, the evaluation map \eqref{eval}
is injective
$$ 0 \lra \cO_C \otimes H^0(C, \Hom(E,F)) \map{ev} \Hom(E,F) \lra
\CC_{\Delta(E)} \lra 0. $$
The cokernel is a sky-scraper sheaf supported at a divisor $\Delta(E)$.
Since $\det \Hom(E,F) = \omega^2$, we have $\Delta(E) \in |\omega^2|$.
This shows that set-theoretically we have $C \cap \Hess(E) = \Delta(E)$.
Let us determine the local equation of $\Hess(E)$ at a point $p \in
\Delta(E)$. We denote by $m$ the multiplicity of $\Delta(E)$ at the 
point $p$. Then, since there is no section of $\Hom(E,F)$ vanishing
twice at $p$ (stability of $E$ and $F$), we have 
$\dim H^0(\Hom(E,F)(-p)) = m$. We choose a basis $\phi_1, \ldots, \phi_m$
of sections of the subspace $H^0(\Hom(E,F)(-p)) \subset H^0(\Hom(E,F))$ and
complete it (if necessary) by $\phi_{m+1}, \ldots,\phi_4$. Let $z$ be a 
local coordinate in an analytic neighbourhood centered at the point $p$.
With these  notations the quadrics $Q_z$ of the net can be written
$$Q_z(\lambda_1,\ldots,\lambda_4) = \det \left( \sum_{i=1}^{4} 
\lambda_i \phi_i(z) \right), $$
where the $\phi_i(z)$ are a basis of the fibre  $\Hom(E,F)_{z}$ for
$z \not= 0$. By construction we have for $1 \leq i \leq m$, $\phi_i(z) = z 
\psi_i(z)$, and the local equation of $\Hess(E)$ is the determinant of the
symmetric $4 \times 4$ matrix
$$ \Hess(E)(z) = \det \left[ B(\phi_i(z), \phi_j(z)) 
\right]_{1\leq i,j \leq 4},$$
where $B$ is the polarization of the determinant. We obtain that 
$\Hess(E)(z)$ is of the form $z^{2m} R(z)$. Hence
$\mult_p(\Hess(E)) \geq 2m$, proving the statement.  
\end{proof}

\begin{definition} \label{disc}
We call the divisor $\Delta(E)$ the discriminant divisor of $E$ and the
rational map $\Delta: \M_0 \lra |\omega^2|$ the discriminant map.
\end{definition}

In the sequel of this paper we will show that the bundle $E$ and
its Hessian curve $\Hess(E)$ are in bijective correspondence (modulo
some discrete structure, which will be defined in section 4.5.2). A first
property is the following: given $E \in \M_0^{reg}$, we associate to the $28$ 
degree two effective divisors $D_{ij}$ \eqref{gbt} on the curve $C$ 
their corresponding secant lines $\overline{D}_{ij} \subset |\omega|^*$.

\begin{prop} \label{secbi}
The secant line $\overline{D}_{ij}$ to the curve $C$ coincides with the
bitangent $\Lambda_{ij}$ to the smooth quartic curve $\Hess(E)$.
\end{prop}

\begin{proof}
Since the bitangent $\Lambda_{ij}$ to $\Hess(E)$ corresponds to the
pencil of quadrics in $\Pi$ containing the line $\overline{x_i x_j}$,
it will be enough to show that $Q_a$ and $Q_b$ belong to $\Lambda_{ij}$,
for $D_{ij} = a+b$, with $a,b \in C$. Consider the vector bundle
map $\pi_i \oplus \pi_j : E \lra L_i \oplus L_j$, where $\pi_i$ and
$\pi_j$ are the natural projection maps. Since $L_i L_j = \cO(D_{ij})$,
the map $\pi_i \oplus \pi_j$ has cokernel $\CC_a \oplus \CC_b$, which
is equivalent to saying that the two linear forms $\pi_{i,a}: E_a \ra
L_{i,a}$ and $\pi_{j,a}: E_a \ra L_{j,a}$ are proportional (same
for $b$). This implies that any map $\phi \in \overline{x_i x_j}$
factorizes at the point $a$ through $\pi_{i,a} = \pi_{j,a}$,
hence $\det \phi_a = 0$. This means that $\overline{x_i x_j}
\subset Q_a$,i.e., $Q_a \in \Lambda_{ij}$ (same for $b$). 
\end{proof}

\subsection{Moduli of $\pp \SL_2$-bundles and the discriminant map $\Delta$}

The finite group $JC[2]$ of $2$-torsion points of $JC$ acts by tensor
product on $\M_0$ and $\M_\omega$. Since Coble's quartic is
Heisenberg-invariant, it is easily seen that the polar map
$\DD : \M_0 \lra \M_\omega$ is $JC[2]$-equivariant,i.e., 
$\DD(E \otimes \alpha) = \DD(E) \otimes \alpha$, $\forall 
\alpha \in JC[2]$. This implies that the constructions we made 
in sections 4.2 and 4.3, namely the projective space $\pp^3 = 
\pp H^0(\Hom(E,F))$, the net of quadrics $\Pi$, its Hessian
curve $\Hess(E)$ and discriminant divisor $\Delta(E)$, only
depend on the class of $E$ modulo $JC[2]$, which we denote 
by $\overline{E}$. It is therefore useful to introduce the
quotient $\N = \M_0/JC[2]$, which can be identified with the
moduli space of semi-stable $\pp \SL_2$-vector bundles with fixed
trivial determinant.  We observe that $\N$ is canonically isomorphic
to the quotient $\M_\omega/JC[2]$. Therefore the $JC[2]$-invariant
polar map $\DD$ descends to a birational involution
\begin{equation} \label{inv}
\overline{\DD}: \N \lra \N. 
\end{equation}
We recall \cite{bls} that the generator $\overline{\LL}$ of 
$\pic(\N) = \ZZ$ pulls-back under the quotient map $q: \M_0 \lra
\N$ to $q^* \overline{\LL} = \LL^4$ and that global sections
$H^0(\N,\overline{\LL}^k)$ correspond to $JC[2]$-invariant 
sections of $H^0(\M_0,\LL^{4k})$. 

\bigskip
\noindent
The Kummer variety $\K_0$ is contained in the singular locus of
$\N$: since the composite map $JC \map{i} \M_0 \map{q} \N$, with $i(L) =
[L \oplus L^{-1}]$, is $JC[2]$-invariant, it factorizes $JC \map{[2]}
JC \map{\overline{i}} \N$, and the image $\overline{i}(JC) \cong
\K_0 \subset \N$.

\bigskip
\noindent
We also recall from \cite{op2} that we have a morphism 
$$ \Gamma: \N \lra |3\Theta|_+ = \pp^{13}, \qquad \overline{E}  \lms
\Gamma(\overline{E}) = \{ L \in \pic^2(C) \ | \ \dim H^0(C, \sym^2(E) 
\otimes L) > 0 \}, $$
which is well-defined since $\Gamma(\overline{E})$ only depends on
$\overline{E}$. The subscript $+$ denotes invariant (w.r.t. $\xi \lms
\omega \xi^{-1}$) theta-functions. When restricted to $\K_0$ the
morphism $\Gamma$ is the Kummer map,i.e, we have a commutative
diagram
$$
\begin{CD}
\K_0 @>Kum>> |2\Theta| = \pp^7 \\
@VVV @VV+\Theta V \\
\N  @>\Gamma >> |3\Theta|_+ = \pp^{13}  
\end{CD}
$$
The main result of \cite{op2} is
\begin{prop}
The morphism $\Gamma: \N \lra |3\Theta|_+$ is given by the
complete linear system $|\overline{\LL}|$,i.e., there
exists an isomorphism $|\overline{\LL}|^* \cong |3\Theta|_+$.
\end{prop}

\begin{rem}
Using the same methods as in \cite{narram}, one can show that $\Gamma:
\N \lra |3\Theta|_+$ is an embedding. We do not use that result.
\end{rem}

Since the open subset $\M_0^{reg}$ is $JC[2]$-invariant, we
obtain that $\M_0^{reg} = q^{-1}(\N^{reg})$. By passing to the
quotient $\N$,  the Aronhold bundles \eqref{Aronholdbun}
determine $36 \cdot 8 = 288$ distinct points $A(\theta,x) :=
\overline{A(\theta,L,x)} \in \N^{reg}$, the exceptional 
bundles \eqref{distbun} determine one point in $\N$, 
denoted by $A_0$, and we obtain a (rational) 
discriminant map \eqref{disc}
$$ \Delta: \N \lra |\omega^2| $$
defined on the open subset $\N^{reg} \setminus \{ A(\theta,x) \}$. We also
note that the $28$ line bundles $L_i L_j$ for $L_i \in \MM_0(E)$ only
depend on $\overline{E}$. 

\begin{rem} \label{remaronhold}
The $288$ points $A(\theta,x)$ are in 1-to-1 correspondence with unordered
Aronhold sets (see \cite{do} page 167),i.e., sets of seven odd 
theta-characteristics $\theta_i$ ($1\leq i \leq 7$) such that
$\theta_i + \theta_j - \theta_k$ is even $\forall i,j,k$. The seven
$\theta_i$ are cut out on the Steinerian curve by the seven lines
$\overline{xx_i}$, where $x,x_i$ are the base points of $\Pi$.  
\end{rem}

\bigskip
The main result of this section is

\begin{prop} \label{proj}
We have a canonical isomorphism $|3\Theta_{|\Theta}|_+ \cong
|\omega^2|$, which makes the right diagram commute
$$
\begin{array}{ccccc}
\K_0 &\subset & \N & \map{\Delta} &  |\omega^2| \\
 & & & & \\
\cap & & \downarrow^{\Gamma} & & \downarrow^{\cong} \\
 & & & & \\
|2\Theta| &  \map{+\Theta} &  |3\Theta|_+ & \map{res_{\Theta}} & 
|3\Theta_{|\Theta}|_+. 
\end{array}
$$
In other words, considering $\N$ (via $\Gamma$) as a subvariety
in $|3\Theta|_+$, the discriminant map $\Delta$ identifies with the
projection with center $|2\Theta| = \mathrm{Span}(\K_0)$, or
equivalently, with the restriction map of $|3\Theta|_+$ to the
Theta divisor $\Theta \subset \pic^2(C)$.
\end{prop} 

\begin{proof}
First we show that the discriminant map $\Delta$ is given by a linear
subsystem of $|\overline{\LL}|$ ($\cong |3\Theta|^*_+$). Consider 
a line bundle $L \in \pic^1(C)$ and the composite map
$$ \psi_L: \pp^3:= \pp_0(L) \lra \M_0 \map{q} \N \map{\Delta} |\omega^2|.$$
Then it will be enough to show that $\psi_L^*(H) \in |\cO_{\pp^3}(4)|$
(since $q^* \overline{\LL} = \LL^4$) for a hyperplane $H$ in $|\omega^2|$.
We denote by $p$ (resp.$q$) the projection of $\pp^3 \times C$ onto
$C$ (resp.$\pp^3$). There exists a universal extension bundle
$\EE$ over $\pp^3 \times C$
\begin{equation} \label{fame} 
0 \lra p^*L^{-1} \lra \EE \lra p^*L \otimes q^*\cO_{\pp^3}(-1) \lra 0
\end{equation}
such that $\forall e \in \pp_0(L)$ the vector bundle $\EE_{|\{e\} \times C}$
corresponds to the extension class $e$. We denote by $\WW \hookrightarrow
\cO_{\pp^3} \otimes H^0(\omega L^2)$ the universal rank $3$ subbundle over
$\pp^3$ and we define the family $\FF$ over $U \times C$ by the
exact sequence
\begin{equation} \label{famf}
0 \lra \left( \FF \otimes p^* L \right)^* \lra q^* \WW \map{ev}
p^*(\omega L^2) \lra 0, 
\end{equation}
where $U$ is the open subset $\pp^3 \setminus C$. We have $\FF_{|\{e\}
\times C} \cong F_e$ (see \eqref{deffe}). Note that $\pic(U) = \pic(\pp^3)$.
It follows immediately from \eqref{fame} and \eqref{famf} that
$\det \EE = q^* \cO(-1)$, $\det \FF = q^* \cO(1) \otimes p^* \omega$,
and that $\det \EE \otimes \FF = p^* \omega^2$. After removing 
(if necessary) the point $A_0$ from $U$ (see Remark 3.4(2)), we
obtain that $\forall e \in U$, $\dim H^0(C, \EE \otimes \FF_{|\{e\}
\times C}) = 4$, hence by the base change theorems, the direct 
image sheaves $q_* (\EE \otimes \FF)$ and $R^1q_* (\EE \otimes \FF)$
are locally free over $U$. Suppose that the hyperplane $H$ consists
of divisors in $|\omega^2|$ containing a point $p \in C$. Then 
$\psi_L^*(H)$ is given by the determinant of the evaluation map 
over $U$ (see \eqref{eval})
$$ q_*(\EE \otimes \FF) \map{ev} \EE \otimes \FF_{|U \times \{p\}}. $$
Since $\det (\EE \otimes \FF_{|U \times \{p\}}) = \cO_U$, 
the result will follow from the equality $\det q_*(\EE \otimes \FF)
= \cO_U(-4)$, which we prove by using some properties of the
determinant line bundles \cite{km}. Given any family of bundles
$\F$ over $U \times C$, we denote the determinant line bundle
associated to the family $\F$ by $\det Rq_* (\F)$. First we observe
that by relative duality \cite{kl} we have 
$$ q_* (\EE \otimes \FF) \map{\sim} \left(R^1q_* (\EE \otimes \FF) 
\right)^*, $$
hence $\det Rq_*(\EE \otimes \FF) = \left( \det q_*(\EE \otimes \FF) 
\right)^{\otimes 2}$. Next we tensor \eqref{fame} with $\FF$
$$ 0 \lra \FF \otimes p^*L^{-1} \lra \EE \otimes \FF \lra 
\FF \otimes p^* L \otimes q^* \cO(-1) \lra 0.$$
Since $\det Rq_*$ is multiplicative, we obtain 
$$ \det Rq_*( \EE \otimes \FF) \cong \det Rq_*(\FF \otimes p^*L^{-1})
\otimes \det Rq_* (\FF \otimes p^*L \otimes q^*\cO(-1)).$$
Again by relative duality we have $\det Rq_*(\FF \otimes p^*L^{-1})
\cong \det Rq_* (\FF \otimes p^*L \otimes q^*\cO(-1))$, hence
(since $\pic(U) = \ZZ$) we can divide by $2$ to obtain
$$ \det q_*( \EE \otimes \FF) \cong \det Rq_* 
(\FF \otimes p^*L \otimes q^*\cO(-1)) \cong \det Rq_*(\FF \otimes p^*L)
\otimes \cO(-2).$$
The last equation holds since $\chi(F_eL) = 2$. Finally, we apply the
functor $\det Rq_*$ to the dual of \eqref{famf}
\begin{eqnarray*}
\det Rq_*(\FF \otimes p^*L) & \cong & \det Rq_* (q^* \WW^*) \otimes
\det Rq_*(p^* \omega L^2)^{-1} \\
 & \cong & \left( \det \WW^* \right)^{\otimes \chi(\cO)} \cong \cO(-2).
\end{eqnarray*}
which proves $\det q_* (\EE \otimes \FF) = \cO(-4)$.

\bigskip

We also deduce from this construction that the exceptional locus
of the rational discriminant map $\Delta$ is the union of the Kummer
variety $\K_0$ , the exceptional bundle $A_0$, and the 
$288$ Aronhold bundles $A(\theta,x)$. Therefore the map $\Delta$ is given
by the composite of $\Gamma$ with a projection map $\pi: 
|\overline{\LL}|^* \cong |3\Theta|_+ \lra |\omega^2|$, whose
center of projection $\ker \pi$ contains $\mathrm{Span}(\K_0) =
|2\Theta|$. In order to show that $\ker \pi = |2\Theta|$, it suffices
(for dimensional reasons) to show that $\Delta$ is dominant:

Consider a general divisor $\delta = a_1 + \cdots + a_8 \in |\omega^2|$
and choose $M \in \pic^2(C)$ such that $a_1 + \cdots + a_4 \in |M^2|$,
or equivalently $a_5 + \cdots + a_8 \in |\omega^2 M^{-2}|$. Using Lemma
\ref{trisecant} we can find a stable $E \in \T_0$ such that 
$[\DD(E)] = [M \oplus \omega M^{-1}]$. We easily deduce from 
Remark 3.4(2) that $\Delta(E) = \delta$.

\bigskip

Finally, we deduce from the natural exact sequence associated to
the divisor $\Theta \subset \pic^2(C)$ 
$$ 0 \lra H^0(\pic^2(C), 2\Theta) \map{+\Theta} H^0(\pic^2(C), 3\Theta)_+
\map{res_\Theta} H^0(\Theta, 3\Theta_{|\Theta})_+ \lra 0 $$
that the projectivized restriction map $res_\Theta$ identifies with
the projection $\pi$.
\end{proof}

\begin{rem}
Geometrically the assertion on the exceptional locus of $\Delta$ given
in the proof means that (we map $\N$ via $\Gamma$ into $|3\Theta|_+$)
$$ \N \cap |2\Theta| = \K_0 \cup \{ A_0 \} \cup \{ A(\theta,x) \},$$
or equivalently, that the $3\theta$-divisors $\Gamma(A_0)$ and
$\Gamma(A(\theta,x))$ are reducible of the form
$$ \Gamma(A_0) = \Theta + \Gamma^{res}(A_0), \qquad
 \Gamma(A(\theta,x)) = \Theta + \Gamma^{res}(A(\theta,x)), $$
where the residual divisors $\Gamma^{res}(A_0)$ and 
$\Gamma^{res}(A(\theta,x))$ lie in $|2\Theta|$. This can be
checked directly.
\begin{itemize}
\item[(1)] exceptional bundle $A_0$: since $\Theta \cong \sym^2 C$,
the inclusion $\Theta \subset \Gamma(A_0)$ is equivalent
to $\dim H^0(C,\sym^2(A_0) \otimes \omega^{-1}(p+q)) >0$, $\forall p,q
\in C$ (here we take $A_0 \in \M_\omega$ see \eqref{distbun}),
or $\dim H^0(C,\sym^2(A_0) (-u-v)) >0$, $\forall u,v \in C$. But this
follows immediately from $\dim H^0(C,A_0) = 3$, which implies that
$\forall u$ there exists a nonzero section $s_u \in H^0(C,A_0(-u))$.
Taking the symmetric product, we obtain $s_u \cdot s_v \in 
H^0(C, \sym^2(A_0)(-u-v))$.

\item[(2)] Aronhold bundles $A(\theta,x)$: similarly we have to
show that $\dim H^0(C,\sym^2(A) \otimes \omega(-p-q)) >0$,
$\forall p,q \in C$ (take $A = A(\theta,L,x) \in \M_0$).
Since $\MM_0(A)$ is invariant under the involution $L_i \mapsto
\theta L_i^{-1}$, we have $\DD(A) = A \otimes \theta$ and
$\dim H^0(C,A \otimes A \otimes \theta) = \dim H^0(C, \sym^2(A) \otimes
\theta) = 4$. Hence $\forall p$ there exists a nonzero section
$s_p \in H^0(C,\End_0(A) \otimes \theta(-p))$ (note that
$\End_0(A) = \sym^2(A)$) and by taking the $\End_0$-part of the
composite section $s_p \circ s_q$, we obtain a nonzero element
of $H^0(C,\sym^2(A) \otimes \omega(-p-q))$.
\end{itemize}

Moreover it can be shown by standard methods that $\sym^2(A_0)$ and
$\sym^2(A(\theta,x))$ are stable bundles. It would be interesting 
to describe explicitly the $2\theta$-divisors $\Gamma^{res}(A_0)$
and $\Gamma^{res}(A(\theta,x))$, which, we suspect, do not lie on the
Coble quartic $\M_0$.
\end{rem}

\subsection{The action of the Weyl group $W(E_7)$}

The aim of this section is to show that the Hessian map (section 4.3),
which associates to a $\pp \SL_2$-bundle $\overline{E} \in \N^{reg}$
the isomorphism class of the smooth curve $\Hess(\overline{E}) \in \M_3$,
is dominant.

\subsubsection{Some group theory related to genus $3$ curves}

We recall here (see e.g. \cite{a}, \cite{do}, \cite{m}) 
the main results on root lattices
and Weyl groups. Let $\Gamma \subset \pp^2$ be a smooth plane quartic
and let $V$ be its associated degree $2$ del Pezzo surface,i.e., the
degree $2$ cover $\pi: V \ra \pp^2$ branched along the curve $\Gamma$. 
We choose an isomorphism (called geometric marking of $V$) of
the Picard group $\pic(V)$, 
\begin{equation} \label{geomarking}
\varphi: \pic(V) \map{\sim} H_7 = \bigoplus_{i=0}^{7} \ZZ e_i, 
\end{equation}
with the hyperbolic lattice $H_7$, such that $\varphi$ is orthogonal
for the intersection form on $\pic(V)$ and the quadratic form on $H_7$
defined by $e_0^2 = 1$; $e_i^2 = -1, \ (i \not=0)$; $e_i \cdot e_j =0, \
(i \not= j)$. The anti-canonical class $-k$ of $V$ equals $3e_0 -
\sum_{i=1}^7 e_i$. We put $e_8 := \sum_{i=1}^7 e_i -2e_0 = e_0 +k$. Then
the $63$ positive roots of $H_7$ are of two types
\begin{equation} \label{roots}
(1) \ \alpha_{ij}= e_i - e_j, \ 1\leq i <j \leq 8, \qquad
(2) \ \alpha_{ijk}= e_0 - e_i - e_j - e_k,  \ 
1\leq i<j<k \leq 7.
\end{equation}
The $28$ roots of type (1) correspond to the $28$ positive roots of the
Lie algebra $\mathfrak{sl}_8$ seen as a subalgebra of the
exceptional Lie algebra $\mathfrak{e}_7$. Similarly the $56$
exceptional lines of $H_7$ are of two types 
\begin{equation} \label{exclines}
(1) \ l_{ij} = e_i + e_j -e_8, \qquad (2) \ l'_{ij} = e_0 - e_i - e_j,
\qquad  1 \leq i <j \leq 8.
\end{equation}
The Weyl group $W(\SL_8)$ equals the symmetric group $\Sigma_8$ and
is generated by the reflections $s_{ij}$ associated to the roots
$\alpha_{ij}$ of type (1). The action of the reflection $s_{ij}$ on the
exceptional lines $l_{pq}$ and $l'_{pq}$ is given by applying 
the transposition $(ij)$ to the indices $pq$. The Weyl group $W(E_7)$
is generated by the reflections $s_{ij}$ and $s_{ijk}$ (associated to
$\alpha_{ijk}$)  and the reflection $s_{ijk}$ acts on the 
exceptional lines as follows
\begin{itemize}
\item  if $|\{i,j,k,8 \} \cap \{ p,q \}|=1$, then $s_{ijk}(l_{pq}) =
l_{pq}$,
\item  if $|\{i,j,k,8 \} \cap \{p,q \}| = 0 \ \text{or} \ 2$, then
$s_{ijk}(l_{pq}) = l'_{st}$ such that $\{p,q,s,t\}$ equals $\{i,j,k,8\}$ or its
complement in $\{1,\ldots,8\}$.
\end{itemize}
Let us consider the restriction map $\pic(V) \map{res} \pic(\Gamma)$
to the ramification divisor $\Gamma \subset V$. Then we have
the beautiful fact (see \cite{do} Lemma 8 page 190) that $res$ 
maps bijectively the $63$ positive roots $\{\alpha_{ij}, \alpha_{ijk}\}$
\eqref{roots} to the $63$ nonzero $2$-torsion points $J\Gamma[2]\setminus
\{0\}$, thus endowing the Jacobian $J\Gamma$ with a level-$2$-structure,i.e.,
a symplectic isomorphism $\psi: J\Gamma[2] \cong \FF_2^3 \times \FF_2^3$
(for the details, see \cite{do} chapter 9). We also observe that 
the partition of $J\Gamma[2]$ into  the two sets $\{ res(\alpha_{ij}) \}$
($28$ points) and $\{ res(\alpha_{ijk}), 0 \}$ ($36$ points) corresponds
to the partition into odd and even points (w.r.t. the level-$2$-structure $\psi$).
Moreover the images of the $56$ exceptional lines \eqref{exclines} are
the $28$ odd theta-characteristics on $\Gamma$, which we denote by
$res(l_{ij}) = res(l'_{ij}) = \theta_{ij}$. Moreover $\pi(l_{ij}) = 
\pi(l'_{ij}) = \Lambda_{ij}$, where $\Lambda_{ij}$ is the bitangent
to $\Gamma$ corresponding to $\theta_{ij}$.

\bigskip
\noindent
Two geometric markings $\varphi, \varphi'$ \eqref{geomarking} differ by an
element $g \in O(H_7) = W(E_7)$ and  their induced level-$2$-structures
$\psi,\psi'$ differ by $\bar{g} \in \Sp(6,\FF_2)$. The restriction
map $W(E_7) \ra \Sp(6,\FF_2)$, $g \mapsto \bar{g}$ is surjective with
kernel $\ZZ/2 = \langle w_0 \rangle = \mathrm{Center}(W(E_7))$. The element
$w_0 \in W(E_7)$ acts as $-1$ on the root lattice, leaves $k$ invariant 
$w_0(k)  = k$ and exchanges the exceptional lines $w_0(l_{ij}) = l'_{ij}$.

\bigskip
\noindent
We also note that $w_0 \notin \Sigma_8 \subset W(E_7)$ and that the injective
composite map $\Sigma_8 \ra W(E_7) \ra \Sp(6,\FF_2)$ identifies $\Sigma_8$
with the stabilizer of an even theta-characteristic.

\subsubsection{Two moduli spaces with $W(E_7)$-action} 

We introduce the $\Sigma_8$-Galois cover $\widetilde{\M_0} \ra \M_0^{reg}$ parametrizing
stable bundles $E \in \M_0^{reg}$ with an order on the $8$ line subbundles 
$\MM_0(E) = \{ L_1,\ldots,L_8 \}$. The group $JC[2]$ acts on $\widetilde{\M_0}$
and we denote the quotient $\widetilde{\M_0}/JC[2]$ by $\widetilde{\N}$,
which is a $\Sigma_8$-Galois cover $\widetilde{\N} \ra \N^{reg}$. The polar map
$\overline{\DD}: \N \ra \N$ \eqref{inv} lifts to a $\Sigma_8$-equivariant
birational involution $\widetilde{\DD}: \widetilde{\N} \ra \widetilde{\N}$.

\bigskip
\noindent
We also consider the moduli space $\P_C$ parametrizing pairs 
$(\Gamma,\varphi)$,
with $\Gamma \subset |\omega|^* = \pp^2$ a smooth plane quartic curve which
satisfies $\Gamma \cap C = 2\Delta$ and $\Delta \in |\omega^2|$, and
$\varphi$ a geometric marking \eqref{geomarking} for the Del Pezzo
surface $V$ associated to $\Gamma$. Then the forgetful map 
$(\Gamma,\varphi) \mapsto \Gamma$ realizes $\P_C$ as a $W(E_7)$-Galois
cover of the space $\R$ of smooth quartic curves  $\Gamma$ satisfying the
above intersection property. Since the general fibre $f^{-1}(\Delta)$ of the 
projection map $\R \map{f} |\omega^2|$ corresponds to the pencil
of curves spanned by the curve $C$ and the double conic $Q^2$ defined by
$Q \cap C = \Delta$, we see that $\R$ is a $\pp^1$-bundle over $|\omega^2|$,
hence rational.

\begin{prop} \label{birat}
The Hessian map (section 4.3) induces a birational map
$$\widetilde{\Hess}: \widetilde{\N} \lra \P_C,$$
which endows $\widetilde{\N}$ with a $W(E_7)$-action. The
action of $w_0$ corresponds to the polar map 
$\widetilde{\DD}$.
\end{prop} 
  
\begin{proof}
Let $\overline{E} \in \widetilde{\N}$ be represented by $E \in \M_0^{reg}$ and
by an ordered set $\MM_0(E) = \{ L_1,\ldots,L_8 \}$. In order to
construct the data $(\Gamma,\varphi)$,  we consider the Del Pezzo
surface $V \map{\pi} \pp^2$ associated to the Hessian curve $\Gamma = \Hess(E)
\subset |\omega|^* =\pp^2$. Since $\Gamma \cap C = 2\Delta(E)$, the
preimage $\pi^{-1}(C) \subset V$ splits into two irreducible
components $C_1 \cup C_2$, with $C_1 = C_2 = C$.  More generally, it can be
shown that the preimage $\pi^{-1}(C \times \R) \subset \mathcal{V}$ has two
irreducible components, where $\mathcal{V} \ra \R$ is the family of Del
Pezzo's parametrized by $\R$. This allows us to choose uniformly 
a component $C_1$. Then by Proposition \ref{secbi} the secant line $\overline{D}_{ij}$
coincides with a bitangent to $\Gamma$. Therefore the preimage 
$\pi^{-1}(\overline{D}_{ij})$ splits into two exceptional lines and we denote
by $l_{ij}$ the line which cuts out the divisor $D_{ij}$ on the curve 
$C_1 = C$. Then the other line $l'_{ij}$ cuts out the divisor $D'_{ij}$ on $C_1$
with $D_{ij} + D'_{ij} \in |\omega|$. Now it is immediate to check that the
classes $e_i = l_{i8}$ for $1 \leq i \leq 7$ and $e_0 = e_i + e_j - l_{ij} -k$ determine
a geometric marking \eqref{geomarking}.

\bigskip
\noindent
Conversely, given $V$ and a geometric marking $\varphi$, we choose a line bundle
$L_8 \in \pic^1(C)$ such that $\omega L_8^2 = e_{0|C=C_1}$. Next we define
$L_i$ for $1 \leq i \leq 7$ by $L_i L_8 = e_{i|C = C_1}$. Then one 
verifies that $l_{ij|C=C_1} = L_iL_j$ and therefore by Remark \ref{reconst}
there exists a bundle $E \in \M_0$ such that $\MM_0(E) = \{L_1,\ldots,L_8\}$.
Since $L_8$ is defined up to $JC[2]$, this construction gives an element of $\widetilde{\N}$.

\bigskip
\noindent
Since the element $\overline{E} \in \widetilde{\N}$ is determined by the $28$
line bundles $L_iL_j$, it will be enough to describe the action of
$\widetilde{\DD}$ and $w_0 \in W(E_7)$ on the $L_iL_j$'s. Suppose
$\widetilde{\DD}(\overline{E}) = \overline{F}$ with $\MM_0(F) =
\{M_1, \ldots,M_8 \}$, then it follows from the equality $\MM_\omega(F) = \MM_0(E)$
(assuming $F = \DD(E)$) that $M_iM_j = \omega L_i^{-1}L_j^{-1}$. On the other hand we have
$w_0(l_{ij}) = l'_{ij}$ and $l_{ij} + l'_{ij} = -k$. Restricting to
$C = C_1$ ($-k_{|C} = \omega$), we obtain that $w_0 = \widetilde{\DD}$.
\end{proof}

\begin{cor} \label{dom}
The morphism $\Hess : \N^{reg} \lra \R$, $\overline{E} \lms \Hess(\overline{E})$
is finite of degree $72$ and $\N^{reg} \lra \M_3$, $\overline{E} \lms
\text{iso class}(\Hess(\overline{E}))$ is dominant.
\end{cor}

\begin{proof}
The first assertion follows from $|W(E_7)/\Sigma_8| = 72$. For the
second it suffices to show that the forgetful map $\R \ra \M_3$ is
dominant: since $\M_3 = |\cO_{\pp^2}(4)|/\pp\GL_3$ a general curve
$C' \in \M_3$ can be represented as a quartic in $|\omega|^* = \pp^2$.
Let $[C],[C'] \in |\cO_{\pp^2}(4)|= \pp^{14}$ denote the quartic equations
of $C$ and $C'$. Projection with center $[C]$ maps $|\cO_{\pp^2}(4)| \lra
|\omega^4|$. We denote by $\mathcal{V}$ the cone with vertex $[C]$ over the
Veronese variety $\mathrm{Ver}|\omega^2| \hookrightarrow |\omega^4|$
and by $\cO$ the closure of the $\pp\GL_3$-orbit of $[C']$. Now since 
$\dim \mathcal{V} = 6$ and $\dim \cO = 8$ ($C'$ general), it follows
that $\mathcal{V}$ and $\cO$ intersect in $\pp^{14}$. By 
construction a point of $\mathcal{V} \cap \cO$ corresponds
to a quartic in $\R$ lying over $C'$.  
\end{proof}

\begin{rem}
The action of the reflection $s_{ijk} \in W(E_7)$ on $\widetilde{\N}$
is easily deduced from its action on the exceptional lines $l_{pq}$ and
$l'_{pq}$ (see section 4.5.1). Representing 
an element $\overline{E} \in \widetilde{\N}$ by $e \in |\omega L_8^2|^*$, it
is easily checked that the restriction of $s_{ijk}$ to $|\omega L_8^2|^*$ is
given by the linear system of quadrics on $|\omega L_8^2|^*$  passing through
the $6$ points $D_{ijk} = D_{i8} + D_{j8} + D_{k8}$. 
In this way we can construct
the $72 = 2 (1+\binom{7}{3})$ bundles in the fibre of $\Hess : \N^{reg} \ra 
\R$.
\end{rem}

\bigskip
\flushleft{Christian Pauly \\
Laboratoire J.-A. Dieudonn\'e \\
Universit\'e de Nice Sophia Antipolis \\
Parc Valrose \\
F-06108 Nice Cedex 02,  France \\
E-mail: pauly@math.unice.fr}

\end{document}